 \documentstyle[graphicx,amscd,amssymb,verbatim,11pt,righttag,color]{amsart}
 \newtheorem{theorem}{Theorem}[section]

 \numberwithin{equation}{section}
 \renewcommand{\rm}{\normalshape}
 
\begin{document}

\title  {Perfect fractal sets with zero Fourier dimension and arbitrarily long arithmetic progressions}

\author{Chun-Kit Lai}

\address{Department of Mathematics, San Francisco State University,
1600 Holloway Avenue, San Francisco, CA 94132.}

 \email{cklai@@sfsu.edu}

\subjclass[2010]{Primary 28A78, 42A38}
\keywords{Arithmetic progression, Fourier dimension, fractals, perfect sets, Moran sets}

 \maketitle

 \begin{abstract}
 By considering a Moran-type construction of fractals on $[0,1]$, we show that for any $0\le s\le 1$, there exists some Moran fractal sets, which is perfect, with Hausdorff dimension $s$ whose Fourier dimension is zero and it contains arbitrarily long arithmetic progressions.
 \end{abstract}
 \section{introduction}

\subsection{Background and Main Results.} Let $E$ be a set inside a locally compact abelian group $G$ (which may be ${\mathbb R}, {\mathbb T}, {\mathbb Z}/n{\mathbb Z}$). It has been widely believed that the decay of the Fourier transform of a measure supported on $E$ and the length of arithmetic progressions (AP) that $E$ can contain is closely related. For example, on ${\mathbb R}$, any sets of positive Lebesgue measure must contain arbitrarily long arithmetic progressions due to the Lebesgue density theorem and it supports measures decaying to zero due to the Riemann-Lebesgue lemma. On the other hand, the uniformity norm of a set $A\subset {\mathbb Z}/n{\mathbb Z}$ is defined as
 $$
 \|A\|_u = \sup_{k\in{\mathbb Z}/n{\mathbb Z}} |\widehat{{\bf 1}_A}(k)|.
 $$
Small enough $\|A\|_u$ will produce 3-term non-trivial AP \cite[Chapter 4 and 10]{TV10}. Variants and generalization of such observation lead to many far reaching consequences, including the Green-Tao Theorem \cite{GT08} on arbitrary long AP on primes.

\medskip
 Throughout the paper, we will focus on the group ${\mathbb R}$.   The Fourier transform of a finite Borel measure $\mu$ on ${\mathbb R}$ is defined to be
 $$
 \widehat{\mu}(\xi) = \int e^{-2\pi i \xi x}d\mu(x).
 $$

 \medskip

    Similar questions on AP and Fourier decay can be asked in the sparse set on ${\mathbb R}$. Fourier decay will be characterized by Fourier dimension. Let $E$ be a Borel set in ${\mathbb R}$, the {\it Fourier dimension} of $E$ is defined to be
 $$
 \mbox{dim}_F(E) = \sup\{\beta: \exists \mu \ \mbox{such that} \ \mbox{spt}\mu = E \ \mbox{and} \  |\widehat{\mu}(\xi)|\le C|\xi|^{-\beta/2}\}.
 $$
 Here, spt$\mu$ is the support of $\mu$. It is well known Fourier dimension is always less than or equal to the Hausdorff dimension \cite[Chapter 4]{Fa}. A set $E$ is called a {\it Salem set} if the Hausdorff dimension and the Fourier dimension of $E$ are equal. It is difficult to construct Salem set deterministically, but it does exist in abundance under many probability models.

 \medskip

 {\L}aba and Pramanik \cite{LP09} initiated the study about the connection between Fourier decay and the existence of AP. They showed, under some additional assumption, the existence of measures with Fourier decay implies that its support contain a 3-term non-trivial AP. Their results generated immense interest investigating how Fourier decay implies the existence of points of some prescribed configuration \cite{C,CIM,S16}.    On the other hand, Shmerkin \cite{S16} recently showed however that there exist Salem sets without any 3-term AP. In this paper, we ask a converse question: {\it Suppose that a perfect set contains an arbitrarily long AP, can we say something about the Fourier dimension of the sets?}

 \medskip

  Recall that a {\it perfect set} on ${\mathbb R}^d$ is a compact set without isolated points. Note that the converse question will be trivial if we consider arbitrary sets. For example, we can take union of a Cantor set with the rational numbers. Then it must contain arbitrarily long AP from the rational numbers, but its Fourier dimension is completely determined by the Cantor sets we are taking. The following is our main conclusion:

 \begin{theorem}
 For any $0\le s\le 1$, there exists a perfect set $E$ whose Hausdorff dimension equals $s$ containing arbitrarily long arithmetic progression and its Fourier dimension equals zero.
\end{theorem}
 \medskip

 The set $E$ is in fact a Cantor-type construction, which is automatically perfect.  More specifically, it is called {\it Moran-type construction} of fractal sets. It means that the contraction ratio at the same level is the same but are allowed to vary at different levels (See Section 2). We show that if at each level, the fundamental intervals we choose are ``biased'' towards one end, then the Moran fractal will have Fourier dimension zero. This is in stark contrast to known results that random Moran construction generically are Salem sets  \cite{S16} (see also \cite[Section 6]{LP09}). By doing so,  we construct perfect sets with arbitrary Hausdorff dimensions that have zero Fourier dimension but it has arbitrary long AP.

 \medskip

\subsection{Open Problems.} This result suggests us some broader open problems in connection to the current research between AP and Fourier decay. Suppose that ${\mathcal E}: = \{E_j\}_{j=1}^{\infty}$ is a countable collection of finite sets on ${\mathbb R}^d$.

 \medskip

 {\bf (Qu1):} {Does there exist any countable collection such that if a perfect set $E$ contains some affine copies of $E_j$, for all $j$, then $E$ supports a measure with decaying Fourier transform?}

 \medskip

 Our theorem showed that $E_j = \{0,1,...,j-1\}$ does not guarantee {\bf (Qu1)} to hold. We may actually ask {\bf (Qu1)} on ${\mathbb R}^1$ in the most general form:

 \medskip

 {\bf (Qu2):} {Suppose that a perfect set $E$ contains affine copies of any finite pattern on ${\mathbb R}^1$, then $E$ supports a measure $\mu$ with decaying Fourier transform? Can the support of $\mu$ be a Salem Set?} 
  
  \medskip

Note that these questions are dealing with sets of zero Lebesgue measure.  In fact, perfect sets of zero Lebesgue measure containing all finite patterns exist and it were first constructed by Erd\"{o}s and Kakutani \cite{EK} in 1957.  Recently, Molter and  Yavicoli \cite{MY}  constructed closed sets of arbitrary Hausdorff dimension on ${\mathbb R}^1$ containing all possible finite patterns.  However, no Fourier analysis on such sets has been studied and {\bf (Qu1)} and {\bf (Qu2)} remains unknown.
 \medskip

  On ${\mathbb R}^2$, one may try
 \begin{equation}\label{eq1}
 E_j = \{(\cos(2\pi k/j,\sin(2\pi k/j)): k=0,1,...,j-1\},
 \end{equation}
 the set of points on ${\mathbb R}^2$ forming vertices of polygons with number of edges equal $j$.

 \medskip

 {\bf (Qu3):} {Suppose that a perfect set $E$ contains some affine copies of $E_j$ in (\ref{eq1}) for all $j$. Does $E$ support a measure $\mu$ with decaying Fourier transform? Is the support of $\mu$ inside $E$ a Salem set?}

 \medskip

\begin{figure}[h]
  \includegraphics[width=8cm]{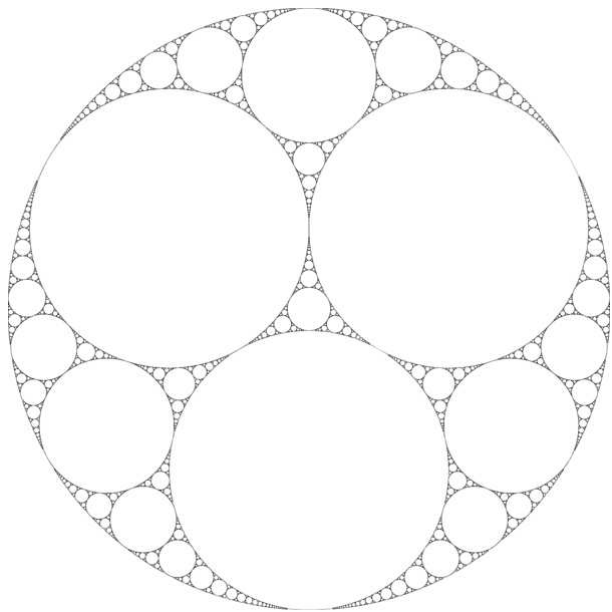}\\
  \caption{Apollonian gasket is formed by deleting the interiors of three circles tangential to the unit circle. Repeating the process by deleting the circle tangential and inside each curvilinear triangle left. The set that remains will be perfect and it contains all $E_j$ in (\ref{eq1}). However, its exact Fourier dimension is largely unknown.}\label{fig}
\end{figure}

  It is clear that the unit circle contains all $E_j$ in (\ref{eq1}), and yet unit circle is a Salem set and it does not give a counterexample to the question. Perfect fractal sets containing all such $E_j$ is abundant in nature. One explicit example is called the {\it Apollonian Gasket} (See \ref{fig}). For a more detailed exposition of Apollonian Gasket, we refer the reader to \cite{GLMWY,K13}. 

\medskip

We also remark that special cases of the Moran measures were considered by the author in the study of fractal measures with exponential bases and frames \cite{GL14,LW16}.

\medskip

 \medskip

\medskip

For the organization of the rest of the paper, we will present our setup and main result in Section 2 and prove our theorems in Section 3.

\medskip

 \section{Setup and Main Results}

 Let $(K_j)_{j=1}^{\infty}$ and $(N_j)_{j=1}^{\infty}$ be sequences of positive integers such that $0<K_j<N_j$ and $N_j\ge 2$. For each $j\ge 1$, let also $B_j$ be subset of $\{0,1,...,N_j-1\}$ satisfying  $\#B_j= K_j$. Moreover, we define
\begin{equation}\label{eqc}
c:= \sup_{j=1,2...}\frac{\max B_j}{N_j}.
 \end{equation}
Define the following multi-index notation:
$$
\Sigma_n = \left\{(b_1,...,b_n): b_j\in B_j\right\}, \ \Sigma^{\ast} = \bigcup_{n=1}^{\infty}\Sigma_n.
$$
For each ${\bf j}: = (b_1,...,b_n)\in \Sigma^{\ast}$, we associate an interval
$$
I_{\bf j} = \left[\sum_{j=1}^{n}\frac{b_j}{N_1...N_j},\sum_{j=1}^{n}\frac{b_j}{N_1...N_j}+\frac{1}{N_1...N_n}\right].
$$
 These intervals form the fundamental intervals of a fractal. The {\it Moran set} associated with $(N_j,B_j)$ is the unique compact set in $[0,1]$ such that
 \begin{equation}\label{eqK}
 E = E(N_j,B_j) = \bigcap_{n=1}^{\infty}\bigcup_{{\bf j}\in \Sigma_n} I_{\bf j}.
 \end{equation}
 $E$ also admits a representation in terms of digit expansion.
 \begin{equation}\label{eqK1}
E = \left\{\sum_{j=1}^{\infty}\frac{b_j}{N_1...N_j}: b_j\in B_j \ \forall j\ge 1\right\}.
 \end{equation}

 Shmerkin \cite[Theorem 2.1]{S16} proved recently that the Moran set in (\ref{eqK1}) is almost surely a Salem set if we choose randomly the digit $B_j$ from $\{0,1,...,N_j-1\}$ with $\#B_j = K_j$ under a natural assumption on $N_j$:
%
  \begin{equation}\label{eqassumption}
  \lim_{j\rightarrow\infty}\frac{\log N_{j+1}}{\log(N_1...N_{j})}=0,
  \end{equation}

       In contrast to the theorem of Shmerkin, we prove if we pick digits biased towards one end, we don't have any decaying Fourier transform.

 \medskip

 \begin{theorem}\label{theorem1}
 Suppose that $c<1/2$ where $c$ is defined in (\ref{eqc}). Then any measures supported on $E$ in (\ref{eqK}) does not decay and thus
 $$
 \mbox{\rm dim}_F(E) = 0.
 $$
 \end{theorem}

The above theorem leads to  the following main theorem.
 \medskip

 \begin{theorem}\label{theorem0}
There exists fractal sets with zero Fourier dimension but it contains arbitrarily long arithmetic progression.
 \end{theorem}

\medskip
The condition $c<1/2$ is sharp in the sense that if all $N_j=2$ and $B_j= \{0,1\}$, then $c = 1/2$ and the measure is the Lebesgue measure on $[0,1]$, whose Fourier transform will decay.

 \medskip

 To avoid such triviality, we consider $N_j\ge 3$ and the consecutive digits set $B_j = \{0,1,...,K_j-1\}$ with $0<K_j<N_j$, so that none of the $E$ generated are intervals. We show that for a larger value of $c$, there exists non-decaying measures on $E$.

\begin{theorem}\label{theorem2}
Let $N_j\ge 3$ and $B_j = \{0,1,...,K_j-1\}$ with $K_j<N_j$. Then

\medskip

\begin{enumerate}
\item if $c<2/3$, {\rm dim}$_F(E)=0$.
\item  if
$$
c<\frac{\sqrt{6}}{\pi}\simeq0.7796968,
$$
 the Fourier transform of the standard measure supported on $E$ which distributes the mass by $\mu(I_{\bf j}) = \frac{1}{K_1...K_n}$ for all ${\bf j}\in \Sigma_n$ and $n\ge 1$ does not decay in the following sense:
 $$
 \limsup_{n\rightarrow \infty} |\widehat{\mu}(n)|>0.
 $$
\end{enumerate}

\end{theorem}

\medskip

We don't know whether dim$_F(E)=0$ for $c$ in between 2/3 and $\sqrt{6}/\pi$. From all the examples we have, to the best of our knowledge, if $E$ supports a measure with non-decaying Fourier transform, then Fourier transform of all other measures supported on $E$ does not decay.

\medskip

Finally, we study the Hausdorff dimension of the fractal sets in Theorem \ref{theorem1} and Theorem \ref{theorem2}. Denote by dim$_H$ the Hausdorff dimension, Using result in \cite{FWW}, we have

\begin{theorem}\label{theorem4}
\begin{enumerate}
\item For the fractal sets in Theorem \ref{theorem1},
$$
\liminf_{j\rightarrow\infty}\frac{\log(K_1...K_j)}{\log(N_1....N_j)+\log(N_{j+1}/K_{j+1})}\le \mbox{\rm dim}_H(E) \le \liminf_{j\rightarrow\infty}\frac{\log(K_1...K_j)}{\log(N_1....N_j)}.
$$
\medskip
\item For the fractal sets in Theorem \ref{theorem2},
$$
 \mbox{\rm dim}_H(E)=\liminf_{j\rightarrow\infty}\frac{\log(K_1...K_j)}{\log(N_1....N_j)+\log(N_{j+1}/K_{j+1})}.
$$
\end{enumerate}
\medskip
In particular, for any $0\le s\le1$, there exists fractal set $E$ with {\rm dim}$_H(E)=s$, {\rm dim}$_F(E)=0$ and it contains arbitrarily long AP.
\end{theorem}

 \medskip

 \section{Proof of the theorems}
  We first prove Theorem \ref{theorem1}. This theorem follows from a classical result of Rajchman. We identify $[0,1]$ as the circle ${\mathbb T}= {\mathbb R}/{\mathbb Z}$. Denote by ${\mathcal M}({\mathbb T})$ the set of all finite Borel measures. Consider the following class of measures in ${\mathcal M}(\mathbb T)$
  $$
  R: = \left\{\mu\in{\mathcal M}({\mathbb T}):   \lim_{|n|\rightarrow\infty} \widehat{\mu}(n)=0\right\}
  $$
  The measure $\mu$ in $R$  are called {\it Rajchman measure}. For more detailed results about Rajchman measures, one can refer to \cite{L95}. Suppose that  ${\mathcal E}$  is a collection of subsets of ${\mathbb T}$, we define also
  $$
  {\mathcal E}^{\perp}: = \{\mu\in{\mathcal M}({\mathbb T}): \forall E\in{\mathcal E}, \ \mu(E)=0\}.
  $$
  For a set $E$ and  an integer $n$,  we define the set
  $$
  nE: = \{nx \ (\mbox{mod} \ 1): x\in E\}.
  $$
  Let also $\overline{E}$ be the closure of $E$ under the standard Euclidean topology. The following theorem  can  be found in \cite{KS64}, \cite{R1922} and \cite{L95}. We provide a proof based on the concept of set of uniqueness for  completeness.

  \begin{theorem}\label{thR}
Let  $$
  {\mathcal H}= \left\{E\subset{\mathbb T}: \exists n_k \uparrow\infty \ \mbox{\rm such that} \ \overline{\bigcup_{k=1}^{\infty}n_k E} \ne {\mathbb T} \right\}.
  $$
  Then $R\subset{\mathcal H}^{\perp}$.  In particular, if $E\in{\mathcal H}$, then dim$_F(E)=0$.
  \end{theorem}

\medskip

\noindent{\it Proof.}  We prove it using the idea of the set of uniqueness. A set $E\subset[0,1]$ is called a {\it set of uniqueness} if
 $$
 \sum_{n=0}^{\infty} (a_n\cos(2\pi nx)+b_n\sin(2\pi n x))=0
 $$
 for all $x\in [0,1]\setminus E$, then $a_n = b_n = 0$.  We know that if $E\in {\mathcal H}$, then $E$ is a set of uniqueness \cite[p.50 Theorem II]{S63}. When $E$ is a set of uniqueness, it does not support any measures whose Fourier transform decays at infinity \cite[p.348 Theorem 6.13]{Z}. Therefore, for any $\mu\in R$, $\mu(E) =0$ for all $E\in {\mathcal H}$. Thus, $R\subset{\mathcal H}^{\perp}$.  Finally, this means also that $E$ does not support any measure with decaying Fourier transform, which shows dim$_F(E)=0$.
\qquad$\Box$
\medskip

 As a simple example of the theorem, for the standard middle-third Cantor set $K_3$, the set $3^nK_3 = K_3$. Hence, $ \overline{\bigcup_{k=1}^{\infty}3^n K_3} = K_3\ne {\mathbb T}$. Thus, dim$_F(K_3)=0$.

 \medskip

 \noindent{\it Proof of Theorem \ref{theorem1}} We show that the Moran set $K$ in (\ref{eqK1}) belongs to ${\mathcal H}$. In fact, we will show that for $n_k = N_1...N_k$,
 \begin{equation}\label{eqclaim}
 \overline{\bigcup_{k=1}^{\infty}n_k E} \subset \left[0,2 c\right].
  \end{equation}
Since $c<1/2$, this union cannot be the whole circle and hence, this will complete the proof by Theorem \ref{thR}.

\medskip

To show that (\ref{eqclaim}) holds, we note from (\ref{eqK1}) that for $n_k = N_1...N_k$
$$
n_k E = \left\{\sum_{j=1}^{\infty} \frac{b_{k+j}}{N_{k+1}...N_{k+j}}: b_{k+j}\in B_{k+j}\right\}.
$$
Using the definition of $c$ and $N_k\ge 2$, we have
\begin{equation}\label{eq3.2+}
\max (n_k E) = \sum_{j=1}^{\infty} \frac{\max B_{k+j}}{N_{k+1}...N_{k+j}} \le  \sum_{j=1}^{\infty} \frac{c}{N_{k+1}...N_{k+j-1}} \le  \sum_{j=1}^{\infty} \frac{c}{2^{j-1}} =2c.
\end{equation}
Hence, $n_k E \subset[0,2c]$ for all $k$. Thus, this shows (\ref{eqclaim}) holds and completes the proof. \qquad$\Box$

\medskip

 \noindent{\it Proof of Theorem \ref{theorem0}}
 We consider $B_j = \{0,1,...,K_j-1\}$, where  $K_j$ satisfies
 \begin{enumerate}
 \item $\lim_{j\rightarrow\infty} K_j=\infty$
 \item $\frac{K_j}{N_j}<\frac12$ for all $j$.
 \end{enumerate}
 Then the second condition and Theorem \ref{theorem1} shows that dim$_F(E(N_j,B_j))=0$. On the other hand, from (\ref{eqK1}), we know that the arithmetic progression
 $$
 \left\{0,\frac{1}{N_1...N_j},...,\frac{K_j-1}{N_1...N_j}\right\}
 $$
 are inside $E(N_j,B_j)$. Condition (1) on $K_j$ implies that we have arbitrarily long arithmetic progression.
\qquad$\Box$

\medskip

We now turn to prove Theorem \ref{theorem2}. The first statement follows from the same proof as Theorem \ref{theorem1}, except in the last step (\ref{eq3.2+}), we use $N_j\ge 3$ to obtain $3/2 c$ in the last step and hence $c<2/3$ will guarantee $\overline{\bigcup_{k=1}^{\infty}n_k E}$ does not cover ${\mathbb T}$.

 \medskip

To prove the second statement, we note that the measure assigning $\mu(I_{\bf j}) = \frac{1}{K_1...K_n}$ is actually an infinite convolution of discrete measures
$$
\mu = \nu_1\ast\nu_2\ast...
$$
where
$$
\nu_n = \frac{1}{K_n}\sum_{j=0}^{K_n-1}\delta_{j/N_1...N_n}
$$
with $\delta_a$ denotes the Dirac measure at $a$.  Suppose that we write
\begin{equation}\label{eq3.1}
\mu_n = \nu_1\ast\nu_2...\ast\nu_n, \ \mu_{>n} = \nu_{n+1}\ast\nu_{n+2}\ast...
\end{equation}
so that $\mu = \mu_n\ast\mu_{>n}$. Note that
\begin{equation}\label{eq3.2}
\widehat{\mu}(N_1...N_n) = \widehat{\mu_n}(N_1...N_n)\widehat{\mu_{>n}}(N_1...N_n) = \widehat{\mu_{>n}}(N_1...N_n)
\end{equation}
because of the integral periodicity that $\widehat{\mu_n}(N_1...N_n)=1$.

\medskip

 \noindent{\it Proof of Theorem \ref{theorem2} (2).} To prove the required statement, we will show that $$\lim_{n\rightarrow\infty}|\widehat{\mu}(N_1...N_n)|>0.$$
  By (\ref{eq3.2}), it suffices to prove that
 $$
 \inf |\widehat{\mu}(N_1...N_n)| =  \inf |\widehat{\mu_{>n}}(N_1...N_n)|>0.
 $$
  We now compute $|\widehat{\mu_{>n}}(\xi)|$ using (\ref{eq3.1}):
  $$
  |\widehat{\mu_{>n}}(\xi)| = \prod_{j=1}^{\infty}|\widehat{\nu_{n+j}}(\xi)|=\prod_{j=1}^{\infty}\left|\frac{1}{K_{n+j}}\sum_{j=0}^{K_{n+j}-1}e^{-2\pi i j\xi/N_1....N_{n+j}}\right|.
    $$
  Hence,
  $$
  \begin{aligned}
  |\widehat{\mu_{>n}}(N_1...N_n)| =& \prod_{j=1}^{\infty}\left|\frac{1}{K_{n+j}}\sum_{j=0}^{K_{n+j}-1}e^{-2\pi i j/N_{n+1}....N_{n+j}}\right|\\
  =& \prod_{j=1}^{\infty}\left|\frac{e^{2\pi i  K_{n+j}/N_{n+1}...N_{n+j}}-1}{K_{n+j}(e^{2\pi i/N_{n+1}...N_{n+j}}-1) }\right|\\
  =& \prod_{j=1}^{\infty}\left|\frac{\sin(\pi  K_{n+j}/N_{n+1}...N_{n+j})}{K_{n+j}\sin(\pi /N_{n+1}...N_{n+j}) }\right|\\
  \ge & \prod_{j=1}^{\infty}\left|\frac{\sin(\pi  K_{n+j}/N_{n+1}...N_{n+j})}{\pi K_{n+j} /N_{n+1}...N_{n+j} }\right| \ (\mbox{by} \ \sin x\le x)\\
  \ge &\prod_{j=1}^{\infty}\left|1-\frac{1}{6}\cdot\left(\frac{\pi  K_{n+j}}{N_{n+1}...N_{n+j}}\right)^2\right| \ (\mbox{by} \ \sin x\ge x -\frac{x^3}{6})
  \end{aligned}
  $$
where the last inequality holds provided that the terms inside are positive, which we are going to check. Note that for $j>1$,
$$
\begin{aligned}
\prod_{j=2}^{\infty}\left|1-\frac{1}{6}\cdot\left(\frac{\pi  K_{n+j}}{N_{n+1}...N_{n+j}}\right)^2\right|\ge& \prod_{j=2}^{\infty}\left|1-\frac{1}{6}\cdot\left(\frac{\pi  }{N_{n+1}...N_{n+j-1}}\right)^2\right| \ (\mbox{by} \ K_{n+1}<N_{n+1})\\
\ge&\prod_{j=2}^{\infty}\left|1-\frac{1}{6}\cdot\left(\frac{\pi  }{3^{j-1}}\right)^2\right|  \ (\mbox{by} \ N_j\ge 3)\\
=&\prod_{j=1}^{\infty}\left(1-\frac{\pi^2}{6\cdot3^{2j}}\right) := c_0
\end{aligned}
$$
and $c_0>0$ because $\sum\frac{\pi^2}{6\cdot3^{2j}}<\infty$. Hence,
$$
  |\widehat{\mu_{>n}}(N_1...N_n)|\ge \left|1-\frac{1}{6}\cdot\left(\frac{\pi  K_{n+1}}{N_{n+1}}\right)^2\right|\cdot c_0\ge \left(1-\frac{c^2\pi^2}{6}\right)\cdot c_0
$$
using the fact that $c= \sup\frac{K_{n+1}}{N_{n+1}}$. Hence, if $c<\sqrt{6}/\pi$, the term above will be positive. This shows that $$
\inf  |\widehat{\mu_{>n}}(N_1...N_n)|>0
$$
completing the proof. \qquad$\Box$

\medskip

We now prove Theorem \ref{theorem4} about the Hausdorff dimension. Under our setting in Section 2, we say that $E$ is a {\it homogeneous Moran set} if for any ${\bf j}\in \Sigma_n$, if we let  $I_{{\bf j}b_1},...,I_{{\bf j}b_{n+1}}$ be the intervals inside $I_{\bf j}$ enumerated from left to right, then the leftmost endpoint of $I_{{\bf j}b_1}$ is the same as $I_{\bf j}$, the rightmost endpoint of $I_{{\bf j}b_{n+1}}$ is the same as $I_{\bf j}$ and the gap between each consecutive intervals are the same.  $E$ is a {\it partial homogeneous Moran set} if for any ${\bf j}\in \Sigma_n$,  then the leftmost endpoint of $I_{{\bf j}b_1}$ is the same as $I_{\bf j}$, and the gap between each consecutive intervals are zero (i.e. intervals are packed towards one end).

\medskip

Feng, Wen and Wu proved that

\begin{theorem}\cite[Theorem 2.1 and Lemma 2.2]{FWW}\label{thFWW}
Given two sequences of numbers, $N_j$ and $K_j$, $j=1,2,...$. The Hausdorff dimension of the homogenous Moran set and the partial homogeneous Moran set are $s_1$ and $s_2$ respectively, where
$$
s_1=\liminf_{j\rightarrow\infty}\frac{\log(K_1...K_j)}{\log(N_1....N_j)}, \ s_2=\liminf_{j\rightarrow\infty}\frac{\log(K_1...K_j)}{\log(N_1....N_j)+\log(N_{j+1}/K_{j+1})}.
$$
Moreover, any other Moran sets $E$ with $\#B_j = K_j$, its Hausdorff dimension must satisfy
$$
s_2\le \mbox{\rm dim}_H(E)\le s_1.
$$
\end{theorem}

Note that if we assume that $N_j$ satisfies
\begin{equation}\label{eq3.3.}
\lim_{j\rightarrow\infty}\frac{\log N_{j+1}}{\log(N_1...N_j)} = 0,
\end{equation}
which is the same as the assumption in (\ref{eqassumption}), then we see that $s_1 = s_2$. Assumption (\ref{eq3.3.}) holds if $N_j$ does not grow too fast. For example, $N_j = j$ will do.
\medskip

\noindent{\it Proof of Theorem \ref{theorem4}.} The first statement follows directly from Theorem \ref{thFWW}. If we pick $B_j = \{0,1,...,K_j-1\}$, then the Moran set $E$ is partial homogeneous, so the second statement also follows from Theorem \ref{thFWW}.

\medskip

We now choose $N_j$ and $K_j$ appropriately so that the last statement holds. Consider first $0<s<1$. Take $K_j = \lfloor N_j^{s}\rfloor$, where $\lfloor x\rfloor$ is the largest integers smaller than or equal to $x$.  Take also $B_j = \{0,1,...,K_j-1\}$. Then $\max B_j/N_j\le N_{j}^{s}/N_j = \frac{1}{N_j^{1-s}}<1/2$ as we choose $N_j$ to be strictly increasing. Hence, Theorem \ref{theorem2} holds. We now compute the Hausdorff dimension. Assuming \eqref{eq3.3.} holds,
$$
\mbox{\rm dim}_H(E) = \liminf_{j\rightarrow\infty}\frac{\log(N_1...N_j)^{s}}{\log(N_1....N_j)}=s.
$$
When $s=1$, we take $K_j = \lfloor N_j/3 \rfloor$. Then $\max B_j/N_j\le 1/3$, the same argument works.

\medskip

Suppose that $s=0$. We choose a sequence $s_j>0$ such that
$$
\lim_{j\rightarrow\infty}s_j = 0, \ \mbox{and} \ \lim_{j\rightarrow\infty} N_j^{s_j} = \infty
$$
(For example, take $s_j = (\ln N_j)^{-1}$). Let $K_j =  \lfloor N_j^{s_j}\rfloor$. Then
$$
\mbox{\rm dim}_H(E) \le \liminf_{j\rightarrow\infty}\frac{\log(N_1...N_j)^{s_j}}{\log(N_1....N_j)}=\liminf_{j\rightarrow\infty}s_j = 0.
$$
Moreover, $K_j$ still tends to infinity with $\max B_j/N_j<1/2$ clearly holds. Hence, Theorem \ref{theorem2} gives us the Moran set has zero Fourier dimension and arbitrarily long AP, completing the proof.

\medskip

\noindent{\bf Acknowledgement.} The author would like to thank professor Tuomas Sahlsten for some  insightful discussions about Fourier dimension and arithmetic progression. He would also like thank professor Malabika Pramanik for an insightful discussion during the MSRI Harmonic Analysis: Introductory Workshop in January 2017 and pointing out  that the reference \cite{EK} is closely related to the research. Finally, the author would also like to thank the anonymous referee for his/her useful comments and suggesting the recent preprint \cite{MY}.

 \end{document}